\documentclass[a4paper]{article}
\usepackage{amsmath}
\usepackage{amssymb}
\newcommand{\be}{\begin{equation}}
\newcommand{\ee}{\end{equation}}
\newcommand{\ra}{\rightarrow}

\newcommand{\bea}{\begin{eqnarray}}
\newcommand{\eea}{\end{eqnarray}}
\newcommand{\beas}{\begin{eqnarray*}}
\newcommand{\eeas}{\end{eqnarray*}}

\newcommand{\R}{\mathbb{R}}
\newcommand{\C}{\mathbb{C}}

\newcommand{\nn}{\nonumber}

\newcommand{\wh}{\widehat}

\newcommand{\ol}{\overline}
\newcommand{\pa}{\partial}
\newcommand{\ti}{\times}
\newcommand{\A}{{\cal A}}
\newcommand{\E}{{\cal E}}
\newcommand{\N}{{\cal N}}
\newcommand{\cC}{{\cal C}}

\newcommand{\cU}{{\cal U}}
\newcommand{\cH}{{\cal H}}

\newcommand{\ka}{{\mathbb K}}

\newcommand{\X}{{\cal X}}

\def\zC{\Psi}
\def\wE{\wh{\E}}
\mathchardef\za="710B  %\alpha
\mathchardef\zb="710C  %\beta
\mathchardef\zg="710D  %\gamma
\mathchardef\zd="710E  %\delta
\mathchardef\zve="710F %\epsilon
\mathchardef\zz="7110  %\zeta
\mathchardef\zh="7111  %\eta
\mathchardef\zvy="7112 %\theta
\mathchardef\zi="7113  %\iota
\mathchardef\zk="7114  %\kappa
\mathchardef\zl="7115  %\lambda
\mathchardef\zm="7116  %\mu
\mathchardef\zn="7117  %\nu
\mathchardef\zx="7118  %\xi
\mathchardef\zp="7119  %\pi
\mathchardef\zr="711A  %\rho
\mathchardef\zs="711B  %\sigma
\mathchardef\zt="711C  %\tau
\mathchardef\zu="711D  %\upsilon
\mathchardef\zvf="711E %\phi
\mathchardef\zq="711F  %\chi
\mathchardef\zc="7120  %\psi
\mathchardef\zw="7121  %\omega
\mathchardef\ze="7122  %\varepsilon
\mathchardef\zy="7123  %\vartheta
\mathchardef\zf="7124  %\varomega
\mathchardef\zvr="7125 %\varrho
\mathchardef\zvs="7126 %\varsigma
\mathchardef\zf="7127  %\varphi
\mathchardef\zG="7000  %\Gamma
\mathchardef\zD="7001  %\Delta
\mathchardef\zY="7002  %\Theta
\mathchardef\zL="7003  %\Lambda
\mathchardef\zX="7004  %\Xi
\mathchardef\zP="7005  %\Pi
\mathchardef\zS="7006  %\Sigma
\mathchardef\zU="7007  %\Upsilon
\mathchardef\zF="7008  %\Phi
\mathchardef\zW="700A  %\Omega

\newcommand{\ep}{\hfill$\Box$}
\newcommand{\bp}{\textit{Proof.-} }

\begin{document}

\title{Local Lie algebra determines base manifold
\footnote{Research supported by the Polish Ministry of Scientific
Research and Information Technology under grant No 2 P03A 020
24.}}
\author{Janusz Grabowski}\maketitle

\newtheorem{re}{Remark}
\newtheorem{theo}{Theorem}
\newtheorem{prop}{Proposition}
\newtheorem{lem}{Lemma}
\newtheorem{cor}{Corollary}
\newtheorem{ex}{Example}

\begin{abstract}It is proven that a local Lie algebra in the
sense of A.~A.~Kirillov determines the base manifold up to a
diffeomorphism provided the anchor map is nowhere-vanishing. In
particular, the Lie algebras of nowhere-vanishing Poisson or
Jacobi brackets determine manifolds. This result has been proven
for different types of differentiability: smooth, real-analytic,
and holomorphic.

\bigskip\noindent
\textit{\textbf{MSC 2000:} Primary 17B66; Secondary 17B65, 53D10,
53D17.}

\medskip\noindent
\textit{\textbf{Key words:} smooth manifolds; vector fields;
diffeomorphisms; Lie algebras; vector bundles; Jacobi brackets.}

\end{abstract}

\bigskip
\centerline{\bf Dedicated to Hideki Omori}
\section{Introduction}
The classical result  of  Shanks and Pursell \cite{PS} states that
the Lie algebra ${\cal X} _c(M)$ of all compactly supported smooth
vector  fields  on    a smooth manifold $M$ determines the
manifold $M$, i.e., the Lie algebras ${\cal X}_c(M_1)$ and ${\cal
X}_c(M_2)$ are isomorphic if and only if $M_1$ and $M_2$ are
diffeomorphic. A similar theorem holds for other complete and
transitive Lie algebras of vector fields \cite{KMO1,KMO2} and for
the Lie algebras of all differential and pseudodifferential
operators \cite{DS,GP}.

There is a huge list of papers in which special geometric
situations (hamiltonian, contact, group invariant, foliation
preserving, etc., vector fields) are concerned. Let us mention the
results of Omori \cite{O} (Ch. X) and \cite{O2}([Ch. XII), or
\cite{A,AG,FT,HM,R,JG3}, for which specific tools were developed
in each case. There is however a case when the answer is more or
less complete in the whole generality. These are the Lie algebras
of vector fields which are modules over the corresponding rings of
functions  (we shall call them  {\it modular}). The standard model
of a modular Lie algebra  of vector fields is  the Lie algebra
${\cal X}({\cal F})$  of  all vector fields tangent to a given
(generalized) foliation ${\cal F}$. If Pursell-Shanks-type results
are concerned in this context, let us recall the work of Amemiya
\cite{Am} and our paper \cite{JG}, where the developed algebraic
approach made it possible to consider analytic cases as well. The
method of Shanks and Pursell consists  of the description  of
maximal ideals in the Lie algebra  ${\cal X}_c(M)$ in  terms  of
the points of $M$: maximal ideals are of the form $\widetilde{p}$
\ for $p\in M,$ where $\widetilde{p}$ consists of vector fields
which are flat at $p.$ This method, however, fails in analytic
cases, since analytic vector fields flat at $p$ are zero on the
corresponding component of $M.$ Therefore in \cite{Am, JG} maximal
finite-codimensional subalgebras are used instead of ideals. A
similar approach is used in \cite{GG} for proving that the Lie
algebras associated with Lie algebroids determine base manifolds.

The whole story for modular Lie algebras of vector fields has been
in a sense finished by the brilliant purely algebraic result of
Skryabin \cite{S}, where one associates the associative algebra of
functions with the Lie algebra of vector fields without any
description of the points of the manifold as ideals. This final
result implies in particular that, in the case when modular Lie
algebras of vector fields contain finite families of vector fields
with no common zeros (we say that they are {\it strongly
non--singular}), isomorphisms between them are generated by
isomorphisms  of corresponding algebras  of functions,  i.e., by
diffeomorphisms   of underlying manifolds.

On the other hand, there are many geometrically interesting Lie
algebras of vector fields which are not modular, e.g. the Lie
algebras of hamiltonian vector fields on a Poisson manifold etc.
For such algebras the situation is much more complicated and no
analog of Skryabin method is known in these cases. In \cite{JG2} a
Pursell-Shanks-type result for the Lie algebras associated with
Jacobi structures on a manifold has been announced. The result
suggests that the concept of a Jacobi structure should be
developed for sections of an arbitrary line bundle rather than for
the algebra of functions, i.e., sections of the trivial line
bundle. This is exactly the concept of {\it local Lie algebra} in
the sense of A.~A.~Kirillov \cite{Ki} which we will call also {\it
Jacobi-Kirillov bundle}.

In the present note we complete the Lie algebroid result of
\cite{GG} by proving that the local Lie algebra determines the
base manifold up to a diffeomorphism if only the anchor map is
nowhere-vanishing (Theorem \ref{t5}). The methods, however, are
more complicated (due to the fact that the Lie algebra of
Jacobi-hamiltonian vector fields is not modular) and different
from those in \cite{GG}. A part of these methods is a modification
of what has been sketched in \cite{JG2}. However, the full
generalization of \cite{JG2} for local Lie algebras on arbitrary
line bundles, i.e., the description of isomorphisms of local Lie
algebras, is much more delicate and we postpone it to a separate
paper. Note also that in our approach we admit different
categories of differentiability: smooth, real-analytic, and
holomorphic (Stein manifolds).

\section{Jacobi modules} What we will call {\it Jacobi module} is
an algebraic counterpart of geometric structures which include
{\it Lie algebroids} and {\it Jacobi structures} (or, more
generally, local Lie algebras in the sense of Kirillov \cite{Ki}).
For a short survey one can see \cite{JG5}, where these geometric
structures appeared under the name of {\it Lie QD-algebroids}.

The concept of a {\it Lie algebroid}  (or  its  pure  algebraic
counterpart -- a {\it   Lie pseudoalgebra}) is one of the most
natural concepts in geometry.

\bigskip\noindent
{\bf Definition 1.} Let $R$ be a commutative and unitary ring, and
let $\A$ be a commutative  $R$-algebra. A {\it Lie pseudoalgebra}
over $R$ and $\A$ is an $\A$-module $\E$ together with a bracket
$[\cdot,\cdot] :\E\ti\E\ra\E$ on the module $\E$, and an
$\A$-module  morphism $\za:\E\ra\text{Der}(\A)$ from $\E$ to the
$\A$-module $\text{Der}(\A)$ of derivations of $\A$, called the
{\it anchor} of $\E$, such that \noindent
\begin{description}
\item{(i)} the bracket on $\E$ is $R$-bilinear, alternating, and
satisfies the Jacobi identity:
$$
[[X,Y],Z]=[X,[Y,Z]]-[Y,[X,Z]].
$$
\item{(ii)} For all  $X,Y\in  \E$  and  all $f\in\A$ we have
\be\label{0}
[X,fY]=f[X,Y]+\za(X)(f)Y;
\ee
\item{(iii)}   $\za([X,Y])=[\za(X),\za(Y)]_c$   for   all
$X,Y\in\E$, where $[\cdot,\cdot]_c$ is the commutator bracket on
${\rm Der}(\A)$.
\end{description}
A {\it Lie algebroid} on a vector bundle $E$ over a base manifold
$M$ is a Lie pseudoalgebra on the $(\R,C^\infty(M))$-module
$\E=Sec(E)$ of smooth sections of $E$. Here the  anchor map is
described  by  a vector  bundle morphism $\za:E\ra TM$ which
induces the bracket homomorphism from $(\E, [\cdot,\cdot])$ into
the Lie algebra $(\X(M),[\cdot,\cdot]_{vf})$ of vector fields on
$M$. In this case, as in the case of any faithful $\A$-module
$\E$, i.e., when $fX=0$ for all $X\in\E$ implies $f=0$, the axiom
(iii) is a consequence of (i) and (ii). Of course, we can consider
Lie algebroids in the real-analytic or holomorphic (on complex
holomorphic bundles over Stein manifolds) category as well.

\smallskip
Lie pseudoalgebras appeared first in a paper by Herz  \cite{He}
but one can find similar concepts under more than a dozen  of
names in the   literature   (e.g. Lie modules,  $(R,A)$-Lie
algebras, Lie-Cartan pairs, Lie-Rinehart algebras, differential
algebras, etc.). Lie algebroids were introduced  by Pradines
\cite{Pr} as infinitesimal parts of differentiable groupoids. In
the same year a book by Nelson \cite{Ne} was published, where a
general theory of Lie modules together with a big part of the
corresponding differential calculus can be found. We also refer to
a survey article by Mackenzie \cite{Ma}.

Note that Lie algebroids on a singleton base space are Lie
algebras. Another canonical example is the tangent bundle $TM$
with the canonical bracket $[\cdot,\cdot]_{vf}$ on the space
$\X(M)=Sec(TM)$ of vector fields.

The property (\ref{0}) of the bracket in the $\A$-module $\E$ can
be expressed as the fact that $ad_X=[X,\cdot]$ is a {\it
quasi-derivation} in $\E$, i.e., an $R$-linear operator $D$ in
$\E$ such that $D(fY)=fD(Y)+\widehat{D}(f)Y$ for any $f\in\A$ and
certain derivation $\widehat{D}$ of $\A$ called the {\it anchor}
of $D$. The concept of quasi-derivation can be traced back to
N.~Jacobson \cite{Ja1,Ja2} as a special case of his {\it
pseudo-linear endomorphism}. It has appeared also in \cite{Ne}
under the name of a {\it module derivation} and used to define
linear connections in the algebraic setting. In the geometric
setting, for Lie algebroids, it has been studied in \cite{Mk}, Ch.
III, under the name {\it covariant differential operator}.

Starting with the notion of Lie pseudoalgebra we obtain the notion
of {\it Jacobi module} when we drop the assumption that the anchor
map is $\A$-linear.

\bigskip\noindent
{\bf Definition 2.} Let $R$ be a commutative and unitary ring, and
let $\A$ be a commutative  $R$-algebra. A {\it Jacobi module} over
$(R,\A)$ is an $\A$-module $\E$ together with a bracket
$[\cdot,\cdot] :\E\ti\E\ra\E$ on the module $\E$, and an
$R$-module  morphism $\za:\E\ra\text{Der}(\A)$ from $\E$ to the
$\A$-module $\text{Der}(\A)$ of derivations of $\A$, called the
{\it anchor} of $\E$, such that (i)-(iii) of Definition 1 are
satisfied. Again, for faithful $\E$, the axiom (iii) follows from
(i) and (ii). This concept is in a sense already present in
\cite{He}, although in \cite{He} it has been assumed that $\A$ is
a field. It has been observed in \cite{He} that every Jacobi
module (over a field) of dimension $>1$ is just a Lie
pseudoalgebra.

\bigskip\noindent
{\bf Definition 3.} (cf. \cite{JG5}) A {\it Lie QD-algebroid} is a
Jacobi module structure on the $(\R,C^\infty(M)$-module
$\E=\text{Sec}(E)$ of sections of a vector bundle $E$ over a
manifold $M$.

\bigskip\noindent
The case $\text{rank}(E)=1$ is special by many reasons and it was
originally studied by A.~A.~Kirillov \cite{Ki}. For a trivial
bundle, well-known examples are those given by Poisson or, more
generally, Jacobi brackets (cf. \cite{Li}). In \cite{Ki} such
structures on line bundles are called {\it local Lie algebras} and
in \cite{Mr} -- {\it Jacobi bundles}. We will refer to them also
as to {\it local Lie algebras} or {\it Jacobi-Kirillov bundles} \
and to the corresponding brackets as to {\it Jacobi-Kirillov
brackets}.

\bigskip\noindent
{\bf Definition 4.}  A {\it Jacobi-Kirillov bundle} ({\it local
Lie algebra} in the sense of Kirillov) is a Lie QD-algebroid on a
vector bundle of rank 1. In other words, a {\it Jacobi-Kirillov
bundle} is a Jacobi module structure on the
$(\R,C^\infty(M))$-module $\E$ of sections of a line bundle $E$
over a smooth manifold $M$. The corresponding bracket on $\E$ we
call {\it Jacobi-Kirillov bracket} and the values of the anchor
map $\za:\E\ra\X(M)$ we call {\it Jacobi-hamiltonian vector
fields}.

\bigskip\noindent
It is easy to see (cf. \cite{JG5}) that any Lie QD-algebroid on a
vector bundle of rank $>1$ must be a Lie algebroid. Of course, we
can consider Lie QD-algebroids and Jacobi-Kirillov bundles in
real-analytic or in holomorphic category as well.

\bigskip\noindent
Since quasi-derivations are particular first-order differential
operators in the algebraic sense, it is easy to see that, for a
Jacobi module $\E$ over $(R,\A)$, the anchor map
$\za:\E\ra\text{Der}(A)$ is also a first-order differential
operator, i.e., \be\label{an} \za(fgX)=f\za(gX)+g\za(fX)-fg\za(X)
\ee
for any $f,g\in\A$ and $X\in\E$. Denoting the Jacobi-hamiltonian
vector field $\za(X)$ shortly by $\wh{X}$, we can write for any
$f,g\in\A$ and $X,Y\in\E$,
\bea\label{n1}
[gX,fY]&=&\wh{gX}(f)Y-f\cdot\wh{Y}(g)X+fg[X,Y]\\
&=&g\cdot\wh{X}(f)Y-\wh{fY}(g)X+fg[X,Y],\nn \eea so that for the
map $\zL_X:\A\ti\A\ra\A$ defined by
$\zL_X(g,f):=\wh{gX}(f)-g\wh{X}(f)$ we have \be\label{n2}
\zL_X(g,f)Y=-\zL_{Y}(f,g)X. \ee The above identity implies clearly
that, roughly speaking, $\text{rank}_\A\E=1$ `at points where
$\zL$ is non-vanishing' (cf. \cite{JG5}) and that \be\label{n2a}
\zL_X(g,f)X=-\zL_{X}(f,g)X. \ee The identity (\ref{n2a}) does not
contain much information about $\zL_X$ if there is `much torsion'
in the module $\E$. But if, for example, there is a torsion-free
element in $\E$, say $X_0$, (this is the case of the module of
sections of a vector bundle), then the situation is simpler. In
view of (\ref{n2a}), $\zL_{X_0}$ is skew-symmetric and, in turn,
by (\ref{n2}) every $\zL_X$ is skew-symmetric. Every $\zL_X$ is by
definition a derivation with respect to the second argument, so,
being skew-symmetric, it is a derivation also with respect to the
first argument. Since in view of (\ref{n1}),
$$[gX,fX]=\left(g\wh{X}(f)-f\cdot\wh{X}(g)+\zL_X(g,f)\right)X,
$$
and since $\zL_X$ and $\wh{X}$ respect the annihilator
$\text{Ann}(X)=\{ f\in\A:fX=0\}$, we get easily the following.

\begin{prop} If $\E$ is a Jacobi module over $(R,\A)$,
then, for every $X\in\E$, the map $\zL_X:\A\ti\A\ra\A$ induces a
skew-symmetric bi-derivation of $\A/{\rm Ann}(X)$, the derivation
$\wh{X}$ of $\A$ induces a derivation of $\A/{\rm Ann}(X)$ and the
bracket
\be \{
\ol{f},\ol{g}\}_X=\zL_X(\ol{f},\ol{g})+\ol{f}\cdot\wh{X}(\ol{g})
-\ol{g}\cdot\wh{X}(\ol{f}), \ee
where $\ol{f}$ denotes the class
of $f\in\A$ in $\A/{\rm Ann}(X)$, is a Jacobi bracket on $\A/{\rm
Ann}(X)$ associated with the Jacobi structure $(\zL_X,\wh{X})$.
Moreover, $\A/{\rm Ann}(X)\ni \ol{f}\mapsto fX\in\E$ is a Lie
algebra homomorphism of the bracket $\{\cdot,\cdot\}_X$ into
$[\cdot,\cdot]$.
\end{prop}
For pure algebraic approaches to Jacobi brackets we refer to
\cite{S2,S3,JG0}.
\begin{cor}\label{c1a} If the $\A$-module $\E$ is generated by torsion-free
elements, then for every $X\in\E$, the map $\zL_X:\A\ti\A\ra\A$ is
a skew-symmetric bi-derivation and the bracket bracket
\be\label{n3} \{ {f},{g}\}_X=\zL_X({f},{g})+{f}\cdot\wh{X}({g})
-{g}\cdot\wh{X}({f}), \ee
is a Jacobi bracket on $\A$ associated
with the Jacobi structure $(\zL_X,\wh{X})$. Moreover, $\A\ni
f\mapsto fX\in\E$ is a Lie algebra homomorphism of the bracket
$\{\cdot,\cdot\}_X$ into $[\cdot,\cdot]$.
\end{cor}

\medskip\noindent
For any torsion-free generated Jacobi module, e.g. a module of
sections of a vector bundle, we have additional identities as
shows the following.
\begin{prop}\label{p1} If the $\A$-module $\E$ is generated by torsion-free
elements, then for all $f_1,\dots,f_m\in\A$, $m\ge 2$, and all
$X,Y\in\E$
$$(a)\qquad\qquad (m-1)[FX,Y]=\sum_{i=1}^m[F_iX,f_iY]-[X,FY]$$
and
$$(b)\qquad\qquad
(m-2)[FX,Y]=\sum_{i=1}^{m-1}[F_iX,f_iY]+[F_mY,f_mX],$$ where
$F=\prod_{i=1}^mf_i$, $F_k=\prod_{i\ne k}f_i$.
\end{prop}
\bp (a) We have (cf. (\ref{n1}))
\beas
\sum_{i=1}^m[F_iX,f_iY]&=&\sum_{i=1}^m\left(
F_i\wh{X}(f_i)Y-f_i\wh{Y}(F_i)X+F[X,Y]+\zL_X(F_i,f_i)Y\right)\\
&=&\wh{X}(F)Y-(m-1)\wh{Y}(F)X+mF[X,Y]+\sum_{i\ne
j}F_{ij}\zL_X(f_j,f_i)Y\\ &=&[X,FY]+(m-1)[FX,Y]+ \sum_{i\ne
j}F_{ij}\zL_X(f_j,f_i)Y,
\eeas
where $F_{ij}=\prod_{k\ne i,j}f_k$. The calculations are based on
the Leibniz rule for derivations:
$$\wh{X}(\prod_{i=1}^mf_i)=\sum_{i=1}^mF_i\wh{X}(f_i),$$
etc. Since, due to Corollary \ref{c1a}, $\zL_X$ is skew-symmetric
and $F_{ij}=F_{ji}$, we have
$$\sum_{i\ne j}F_{ij}\zL_X(f_j,f_i)Y=0$$
and (a) follows.

\medskip\noindent (b) In view of (a), we have
\beas
(m-2)[FX,Y]&=&\sum_{i=1}^m[F_iX,f_iY]-[X,FY]-[FX,Y]\\
&=&\sum_{i=1}^{m-1}[F_iX,f_iY]+[F_mX,f_mY]-[X,F_mf_mY]-[F_mf_mX,Y].
\eeas
But
$$[F_mX,f_mY]-[X,F_mf_mY]-[F_mf_mX,Y]=[F_mY,f_mX]$$
is a particular case of (a).

\ep
\section{Useful facts about associative algebras}
In what follows, $\A$ will be an associative commutative unital
algebra over a field $\ka$ of characteristic 0. Our standard model
will be the algebra $\cC(N)$ of class $\cC$ functions on a
manifold $N$ of class $\cC$, $\cC=C^\infty,C^\zw,\cH$. Here
$C^\infty$ refers to the smooth category with $\ka=\R$, $C^\zw$ --
to the $\R$-analytic category with $\ka=\R$, and $\cH$ -- to the
holomorphic category of Stein manifolds with $\ka=\C$ (cf.
\cite{JG,AG}). All manifolds are assumed to be paracompact and
second countable. It is obvious what is meant by a Lie
QD-algebroid or a Jacobi-Kirillov bundle of class $\cC$. The rings
of germs of class $\cC$ functions at a given point are noetherian
in analytic cases that is no longer true in the $C^\infty$ case.
However, all the algebras $\cC(N)$ are in a sense noetherian in
finite codimension. To explain this, let us start with the
following well-known observation.
\begin{theo}\label{t0} Every maximal finite-codimensional ideal of $\cC(N)$
is of the form $\ol{p}=\{ f\in\cC(N):f(p)=0\}$ for a unique $p\in
N$ and $\ol{p}$ is finitely generated.
\end{theo}
\bp The form of such ideals is proven e.g. in \cite{JG},
Proposition 3.5. In view of embedding theorems for all types of
manifolds we consider, there is an embedding
$f=(f_1,\dots,f_n):N\ra\ka^n$, $f_i\in\cC(N)$. Then, the ideal
$\ol{p}$ \ is generated by $\{ f_i-f_i(p)\cdot 1:i=1,\dots,n\}$.
In the smooth case it is obvious, in the analytic cases it can be
proven by means of some coherent analytic sheaves and methods
parallel to those in \cite{JGd}, Note 2.3.

\ep

\bigskip\noindent
{\bf Remark.} Note that in the case of a non-compact $N$ there are
maximal ideals of $\cC(N)$ which are not of the form $\ol{p}$.
They are of course of infinite codimension. It is not known if the
above theorem holds also for manifolds which are not second
countable (cf. \cite{JG6}).

\bigskip\noindent
For a subset $B\subset\A$, by $\text{Sp}(\A,B)$ we denote the set
of those maximal finite-codimensional ideals of $\A$ which contain
$B$. For example, $\text{Sp}(\A,\{ 0\})$ is just the set of all
maximal finite-codimensional ideals which we denote shortly by
$\text{Sp}(\A)$. Put $\overline{B}=\bigcap_{I\in\text{Sp}(A,B)}I$.
For an ideal $I\subset\A$, by $\sqrt{I}$ we denote the radical of
$I$, i.e., $$\sqrt{I}=\{ f\in\A:f^n\in I,\ \text{for some }
n=1,2,\dots\}.$$ The following easy observations will be used in
the sequel.
\begin{theo}\label{t1} \begin{description}

\item{(a)} If $I$ is an ideal of codimension $k$ in $\A$, then
$\sqrt{I}=\overline{I}$ and $(\overline{I})^k\subset I$.

\item{(b)} Every finite-codimensional prime ideal in $\A$ is
maximal.

\item{(c)} If a derivation $D\in\text{Der}(\A)$ preserves a
finite-codimensional ideal $I$ in $\A$, then
$X(\A)\subset\overline{I}$.

\item{(d)} If $I_1,\dots,I_n$ are finite-codimensional and
finitely generated ideals of $\A$, then the ideal
$I_1{\cdot}{\dots}{\cdot}I_n$ is finite-codimensional and finitely
generated.
\end{description}
\end{theo}
\bp (a) The descending series of ideals
$$I+\ol{I}\supset I+(\ol{I})^2\supset\dots$$
stabilizes at $k$th step at most, so
$I+(\ol{I})^k=I+(\ol{I})^{k+1}$. Applying the Nakayama's Lemma to
the finite-dimensional module $(I+(\ol{I})^k)/I$ over the algebra
$A/I$, we get $(I+(\ol{I})^k)/I=\{ 0\}$, i.e., $(\ol{I})^k\subset
I$, thus $\ol{I}\subset\sqrt{I}$. Since for all
$J\in\text{Sp}(A,I)$ we have $\sqrt{I}\subset\sqrt{J}=J$, also
$\ol{I}\supset\sqrt{I}$.

\smallskip\noindent
(b) If $I$ is prime and finite-codimensional, $\sqrt{I}=I$ and
$\sqrt{I}=\ol{I}$ by (a). But a finite intersection of maximal
ideals is prime only if they coincide, so $\ol{I}=J$ for a single
$J\in\text{Sp}(\A)$.

\smallskip\noindent
(c) By Lemma 4.2 of \cite{JG}, $D(I)\subset I$ for a
finite-codimensional ideal $I$ implies $D(\A)\subset J$ for each
$J\in\text{Sp}(\A,I)$.

\smallskip\noindent
(d) It suffices to prove (d) for $n=2$ and to use the induction.
Suppose that $I_1,I_2$ are finite-codimensional and finitely
generated by $\{ u_i\}$ and $\{ v_j\}$, respectively. It is easy
to see that $I_1\cdot I_2$ is generated by $\{ u_i\cdot v_j\}$ and
that $I_1\cdot I_2$ is finite-codimensional in $I_1$. Indeed, if
$c_1,\dots,c_k\in\A$ represent a basis of $A/I_2$, then $\{
c_lu_i\}$ represent a basis of $I_1/(I_1\cdot I_2)$.

\ep
\begin{theo} For an associative commutative unital algebra $\A$
the following are equivalent:
\begin{description}
\item{(a)} Every finite-codimensional ideal of $A$ is finitely
generated.

\item{(b)} Every maximal finite-codimensional ideal of $A$ is
finitely generated.

\item{(c)} Every prime finite-codimensional ideal of $A$ is
finitely generated.
\end{description}
\end{theo}
\bp (a) $\Rightarrow$ (b) is trivial, (b) $\Rightarrow$ (c)
follows from Theorem \ref{t1} (b), and (c) $\Rightarrow$ (a) is a
version of Cohen's Theorem for finite-dimensional ideals.

\ep

\medskip\noindent
{\bf Definition 5.} We call an associative commutative unital
algebra $\A$ {\it noetherian in finite codimension} if one of the
above (a), (b), (c), thus all, is satisfied.

\medskip\noindent
An immediate consequence of Theorem \ref{t0} is the following.
\begin{theo} The algebra $\A=\cC(N)$ is noetherian in finite
codimension.
\end{theo}

\section{Spectra of Jacobi modules}

Let us fix a Jacobi module $(\E,[\cdot,\cdot])$ over $(\ka,\A)$.
Throughout this section we will assume that $\E$ is finitely
generated by torsion-free elements and that $\A$ is a noetherian
algebra in finite codimension over a field $\ka$ of characteristic
0. The $(\ka,\cC(N))$-modules of sections of class $\cC$ vector
bundles over $N$ can serve as standard examples.

For $L\subset\E$, by $\wh{L}$ denote the image of $L$ under the
anchor map: $\wh{L}=\{\za(X):X\in L\}\subset\text{Der}(\A)$. The
set $\wh{\E}$ is a Lie subalgebra in $\text{Der}(\A)$ with the
commutator bracket $[\cdot,\cdot]_c$ and we will refer to
$\wh{\E}$ as to the Lie algebra of `Jacobi-hamiltonian vector
fields'. The main difference with the `modular' case (in
particular, with that of Pursell and Shanks \cite{PS}) is that
$\wh{\E}$ is no longer, in general, an $\A$-module, so we cannot
multiply by `functions' inside $\wh{\E}$. However, we still can
try to translate some properties of the Lie algebra
$(\E,[\cdot,\cdot])$ into the properties of the Lie algebra
$\wh{\E}$ of Jacobi-hamiltonian vector fields by means of the
anchor map and to describe some `Lie objects' in $\E$ or $\wh{\E}$
by means of `associative objects' in $\A$.

\medskip
The {\it spectrum of the Jacobi module} $\E$, denoted by
$\text{Sp}(\E)$, will be the set of such maximal
finite-codimensional Lie subalgebras in $\E$ that not contain
finite-codimensional Lie ideals of $\E$. In nice geometric
situations, $\text{Sp}(\E)$ will be interpreted as a set of points
of the base manifold at which the anchor map does not vanish. Note
that the method developed in \cite{GG} for Lie pseudoalgebras
fails, since the Lemma 1 therein in no longer true for Jacobi
modules. In fact, as easily shows the example of a symplectic
Poisson bracket on a compact manifold, $[\E,\E]$ may include no
non-trivial $\A$-submodules of $\E$. Therefore we will modify the
method from \cite{JG1} where Poisson brackets have been
considered.

Let us fix some notation. For a liner subspace $L$ in $\E$ and for
$J\subset\A$, denote
\begin{itemize}
\item $\N_L=\{ X\in\E:[X,L]\subset L\}$ -- the Lie normalizer of
$L$; \item $U_L=\{ X\in\E:[X,\E]\subset L\}$;
\item $I(L)=\{ f\in\A:\forall X\in\E\ [fX\in L]\}$ -- the
largest associative ideal $I$ in $\A$ such that $I\E\subset L$;
\item  $\E_J=\{ X\in\E:\wh{X}(A)\subset J\}$.
\end{itemize}

\medskip\noindent
It is an easy excercise to prove the following proposition (cf.
\cite{JG1}, Theorem 1.6).
 \begin{prop} If $L$ is a Lie subalgebra in $\E$, then $\N_L$ is a
 Lie subalgebra containing $L$, the set $U_L$ is a Lie ideal in $\N_L$, and
 $\wh{\N_L}(I(U_L))\subset I(U_L)$.
 \end{prop}
 Choose now generators $X_1,\dots,X_n$ of $\E$ over $\A$. For a fixed
 finite-codimensional Lie subalgebra $L$ in $\E$ put
 $U_i=\{ f\in\A: fX_i\in U_L\}$ and $U=\bigcap_{i=1}^nU_i$. Since
 $U_L$ is clearly finite-codimensional in $\E$, all $U_i$ are
 finite-codimensional in $\A$, so is $U$.
 \begin{lem}\label{l1}\

 \begin{description}
 \item{(a)\ } \ $[U^mX_j,X_k]\subset L$ for all $j,k=1,\dots,n$ and $m\ge
 3$.
\item{(b)\ } \ $[U^mX_j,U^lX_k]\subset L$ for all $j,k=1,\dots,n$
and $m,l\ge 1$.
\end{description}
\end{lem}
\bp (a) Take $f_1,\dots,f_m\in U$. Since $f_iX_k\in U_L$,
Proposition \ref{p1} (b) implies $[f_1\cdots f_mX_j,X_k]\in L$.

\smallskip\noindent
(b) The inclusion is trivial for $l=1$, so suppose $l\ge 2$. Take
$f_1,\dots,f_m\in U$, $f_{m+1}\in U^l$ and put $F=f_1\dots
f_{m+1}$, $F_i=\prod_{r\ne i}f_r$. By Proposition \ref{p1} (b)
$$[f_1\cdots f_mX_k,f_{m+1}X_j]=(m-1)[FX_j,X_k]-
\sum_{i=1}^m[F_iX_j,f_iX_k].$$ Since $F\in U^{m+l}$, according to
(a), $[FX_j,X_k]\in L$ and $[F_iX_j,f_iX_k]\subset[\E,U_L]\subset
L$, so the lemma follows.

\ep
\begin{theo}\label{t3} The ideal $I(U_L)$ is finite-codimensional
in $\A$ provided $L$ is a finite-codimensional Lie subalgebra in \
$\E$.
\end{theo}
 \bp Let $\cU$ be the associative subalgebra in $\A$ generated by
 $U$. It is finite-codimensional and, in view of Lemma \ref{l1} (b),
 $[\cU X_j,\cU X_k]\subset L$. Being finite-codimensional in
 $\A$, the associative subalgebra $\cU$ contains a
 finite-codimensional ideal $J$ of $\A$ (cf. \cite{JG1},
 Proposition 2.1 b)). Hence $[JX_j,JX_k]\subset L$ and, since
 $X_i$ are generators of $\E$, $[J\E,J\E]\subset L$. Note that we do
 not exclude the extremal case $\cU=\A=J$. Applying the
 identity
 $$[f_1f_2X,Y]=[f_2X,f_1Y]+[f_1X,f_2Y]-[f_1f_2Y,X]$$
 for $f_1,f_2\in J$, $X\in U_L$, we see that $J^2U_L\subset U_L$.
 In particular, $J^2UX_i\subset U_L$ for all $i=1,\dots,n$, so
 $J^2U\subset U$ and hence $J^2\cU\subset U$ and $J^3\E\subset
 U_L$. Consequently $J^3\subset I(U_L)$. Since $J$ is
 finite-codimensional and finitely generated, $J^3$ is
 finite-codimensional (Theorem \ref{t1} (d)), so $I(U_L)$ is
 finite-codimensional.

 \ep

\medskip\noindent
Denote $\text{Sp}_\E(\A)$ the set of these maximal
finite-codimensional ideals $I\subset\A$ which do not contain
$\wh{\E}(\A)$, i.e., $\E_I\ne \E$. Geometrically,
$\text{Sp}_\E(\A)$ can be interpreted as the support of the anchor
map. Recall that $\text{Sp}(\E)$ is the set of these maximal
finite-codimensional Lie subalgebras in $\E$ which do not contain
finite-codimensional Lie ideals.
\begin{theo} The map $J\mapsto\E_J$ constitutes a bijection of
$\text{Sp}_\E(\A)$ with $\text{Sp}(\E)$. The inverse map is
$L\mapsto\sqrt{I(L)}$.
\end{theo}
\bp Let us take $J\in\text{Sp}_\E(\A)$. In view of (\ref{an}),
$J^2\E\subset\E_J$ which implies that $\E_J$ is
finite-codimensional, as $J^2$ is finite-codimensional and $\E$ is
finitely generated.

We will show that $\E_J$ is maximal. Of course, $\E_J\ne\E$ and
$\E_J$ is of finite codimension, so there is a maximal Lie
subalgebra $L$ containing $\E_J$. We have
$$J^2\E\subset\E_J\subset L\ \Rightarrow\ J^2\subset I(L)\
\Rightarrow\ J\subset\sqrt{I(L)}\ \Rightarrow\ J=\sqrt{I(L)}.
$$
Moreover, $I(L)$ is finite-codimensional, and since, due to
(\ref{0}), $\wh{L}(I(L))\subset I(L)$, then, by Theorem \ref{t1}
(c), $\wh{L}(A)\subset J$, i.e., $L\subset\E_J$ and finally
$L=\E_J$.

Finally, suppose $P$ is a finite-codimensional Lie ideal of $\E$
contained in $\E_J$. Then $U_P$ is a Lie ideal in $\E$ of finite
codimension and, according to Theorem \ref{t3}, $I(U_P)$ is a
finite-codimensional ideal in $\A$. Since $\wh{\E}(I(U_P))\subset
I(U_P)$, and since $I(U_P)\subset I(U_L)\subset J$, we have
$\wh{\E}(A)\subset J$, i.e., $\E=\E_J$; a contradiction.

\medskip
Suppose now that $L\in\text{Sp}(\E)$. Observe first that $\N_L=L$,
since otherwise $L$ would be a Lie ideal, that would, in turn,
imply $U_L\subset L$ and $I(U_L)\subset I(L)$. Since $U_L$ is
finite-codimensional, Theorem \ref{t3} shows that $I(L)$ is
finite-codimensional. Exactly as above we show that
$\wh{L}(\A)\subset\sqrt{I(L)}$, i.e., $L\subset\E_J$, where
$J=\sqrt{I(L)}$. By Theorem \ref{t1} (a), $J^k\subset I(L)$ for
some $k$, so if we had $\E_J=\E$, then $J^k\cdot\E$ would be a
finite-codimensional Lie ideal contained in $L$. Thus $\E_J\ne\E$
and there is $I\in\text{Sp}(\A,J)$ with $\E_I\ne\E$. We know
already that in this case $\E_I$ is maximal. Since
$L\subset\E_J\subset\E_I$ and $L$ is maximal, we have $L=\E_I$ and
$I=J=\sqrt{I(L)}$.

\ep

\begin{cor}\label{c1} Let $(\E,[\cdot,\cdot])$ be a Lie
QD-algebroid of class $\cC$ (i.e., a Jacobi module over
$(\ka,\cC(N))$ of class $\cC$ sections of a class $\cC$ vector
bundle) over a class $\cC$ manifold $N$. Let $S\subset N$ be the
open support of the anchor map, i.e., $S=\{ p\in N: \wh{X}(p)\ne 0
\text{\ for some\ } X\in\E\}$. Then the map $p\mapsto p^*=\{
X\in\E:\wh{X}(p)=0\}$ constitutes a bijection of \ $S$ with
$\text{Sp}(\E)$.
\end{cor}
Let $\wh{\E}$ be the image of the anchor map
$\za:\E\ra\text{Der}(\A)$. By definition of a Jacobi module,
$\wh{\E}$ is a Lie subalgebra in $(\text{Der}(\A),
[\cdot,\cdot]_c)$. Since $\za:\E\ra\wh{\E}$ is a surjective Lie
algebra homomorphism, it induces a bijection of $\text{Sp}(\E)$
onto $\text{Sp}(\wh{\E})$, $L\mapsto\wh{L}=\za(L)$. Thus we get
the following.
\begin{cor}\label{c2} Let $(\E,[\cdot,\cdot])$ be a Lie
QD-algebroid of class $\cC$ over a class $\cC$ manifold $N$. Let
$S\subset N$ be the open support of the anchor map, i.e., $S=\{
p\in N: \wh{X}(p)\ne 0 \text{\ for some\ } X\in\E\}$. Then the map
$p\mapsto \wh{p}=\{ \zx\in\wh{\E}:\zx(p)=0\}$ constitutes a
bijection of \ $S$ with $\text{Sp}(\wh{\E})$.
\end{cor}

\section{Isomorphisms}
It is clear that any isomorphism $\zC:\E_1\ra\E_2$ of the Lie
algebras associated with Jacobi modules $\E_i$ over $(R_i,\A_i)$,
$i=1,2$, induces a bijection
$\zc:\text{Sp}(\E_2)\ra\text{Sp}(\E_1)$. Since the kernels $K_i$
of the anchor maps $\za_i:\E_i\ra\wh{\E_i}$ are the intersections
$$K_i=\bigcap_{L\in\text{Sp}(\E_i)}L,\quad i=1,2,
$$
$\zC(K_1)=K_2$, so $\zC$ induces a well-defined isomorphisms
$$\wh{\zC}:\wh{\E_1}\ra\wh{\E_2}, \quad \wh{\zC}(\wh{X})=\wh{\zC(X)}$$
with the property
\be\label{i1}\wh{L}\in\text{Sp}(\wh{\E_1})\Leftrightarrow
\wh{\zC}(\wh{L})\in\text{Sp}(\wh{\E_2}).
\ee
\begin{prop} If the Lie algebras $(\E_i,[\cdot,\cdot]_i)$,
associated with Jacobi modules $\E_i$, $i=1,2$, are isomorphic,
then the Lie algebras of Jacobi-hamiltonian vector fields \
$\wh{\E_i}$, $i=1,2$, are isomorphic.
\end{prop}
The following theorem describes isomorphisms of the Lie algebras
of Jacobi-hamiltonian vector fields.
\begin{theo}\label{t5} Let $(\E_i,[\cdot,\cdot]_i)$ be a Lie
QD-algebroid of class $\cC$, over a class $\cC$ manifold $N_i$,
and let $S_i\subset N_i$ be the (open) support of the anchor map
$\za_i:\E_i\ra\wh{\E}_i$, $i=1,2$. Then every isomorphism of the
Lie algebras of Jacobi-hamiltonian vector fields
$\zF:\wE_1\ra\wE_2$ is of the form $\zF(\zx)=\zf_*(\zx)$ for a
class $\cC$ diffeomorphism $\zf:S_1\ra S_2$.
\end{theo}
\begin{cor} If the Lie algebras associated with Lie
QD-algebroids $E_i$ of class $\cC$, over class $\cC$ manifolds
$N_i$, $i=1,2$, are isomorphic, then the (open) supports
$S_i\subset N_i$ of the anchor maps $\za_i:\E_i\ra\wh{\E}_i$,
$i=1,2$, are $\cC$-diffeomorphic. In particular, $N_1$ and $N_2$
are $\cC$-diffeomorphic provided the anchors are
nowhere-vanishing.
\end{cor}
{\it Proof of Theorem \ref{t5}.- } According to Corollary
\ref{c2}, the isomorphism $\zF$ induces a bijection $\zf:S_1\ra
S_2$ such that, for every $\zx\in\wE_1$ and every $p\in S_1$,
\be\label{i2} \zx(p)=0\Leftrightarrow\zF(\zx)(\zf(p))=0.
\ee
First, we will show that $\zf$ is a diffeomorphism of class $\cC$.
For, let $f\in\cC(N_1)$. Since the anchor map is a first-order
differential operator, for every $X\in\E_1$ we have
$\wh{f^2X}=2f\cdot\wh{fX}-f^2\cdot\wh{X}$. In particular, for any
$p\in N_1$,
$$\wh{f^2X}(p)-2f(p)\wh{fX}(p)+f^2(p)\wh{X}(p)=0,
$$
so that, due to (\ref{i2}), \be\label{i3}
\zF(\wh{f^2X})(\zf(p))=2f(p)\zF(\wh{fX})(\zf(p))
-f^2(p)\zF(\wh{X})(\zf(p)).
\ee
We can rewrite (\ref{i2}) in the form \be\label{i4}
\zF(\wh{f^2X})=2(f\circ\zc)\cdot\zF(\wh{fX})
-(f\circ\zc)^2\cdot\zF(\wh{X}),
\ee
where $\zc=\zf^{-1}$ and the both sides of (\ref{i4}) are viewed
as vector fields on $S_2$. In a similar way one can get
\be\label{i5}
\zF(\wh{f^3X})=3(f\circ\zc)^2\cdot\zF(\wh{fX})
-2(f\circ\zc)^3\cdot\zF(\wh{X}).
\ee
To show that $f\circ\zc$ is of class $\cC$, choose $q\in S_2$ and
$X\in\E_1$ such that $\zF(\wh{X})(q)\ne 0$. Then we can choose
local coordinates $(x_1,\dots,x_n)$ around $q$ such that
$\zF(\wh{X})=\pa_1=\frac{\pa}{\pa x_1}$. If $a$ is the first
coefficient of the vector field $\zF(\wh{fX})$ in these
coordinates, we get out of (\ref{i4}) and (\ref{i5}) that
$(f\circ\zc)^2-2a(f\circ\zc)$ and $2(f\circ\zc)^3-3a(f\circ\zc)^2$
are of class $\cC$ in a neighbourhood of $q$. But \be\label{i6}
(f\circ\zc)^2-2a(f\circ\zc)=(f\circ\zc-a)^2-a^2
\ee
and
\be\label{i7}
2(f\circ\zc)^3-3a(f\circ\zc)^2=2(f\circ\zc-a)^3+3a(f\circ\zc-a)^2-a^2,
\ee
so $(f\circ\zc-a)^2$ and $(f\circ\zc-a)^3$ are functions of class
$\cC$ in a neighbourhood of $q$, as the function $a$ is of class
$\cC$. Now we will use the following lemma which proves that
$f\circ\zc-a$, thus $f\circ\zc$, is of class $\cC$.
\begin{lem}\label{l2} If $g$ is a $\ka$-valued function in a neighbourhood
of $0\in\ka^n$ such that $g^2$ and $g^3$ are of class $\cC$, then
$g$ is of class $\cC$.
\end{lem}
\bp In the analytic cases the lemma is almost obvious, since
$g=\frac{g^3}{g^2}$ is a meromorphic and continuous function. In
the smooth case the Lemma is non-trivial and proven in \cite{Jo}.

\ep

To finish the proof of the theorem, we observe that $f\circ\zc$ is
of class $\cC$ for all $f\in\cC(N_2)$ implies that $\zc$, thus
$\zf=\zc^{-1}$, is of class $\cC$ and we show that $\zF=\zf_*$ or,
in other words, that $\wh{Y}(f)\circ\zc=\zF(\wh{Y})(f\circ\zc)$
for all $f\in\cC(N_1)$ and all $Y\in\E_1$. Indeed, for arbitrary
$f\in\cC(N_1)$ and $X,Y\in\E_1$, the bracket of vector fields
$[\wh{Y},\wh{f^2X}]$ reads
\beas[\wh{Y},\wh{f^2X}]&=&[\wh{Y},2f\cdot\wh{fX}-f^2\cdot\wh{X}]\\
&=&
2\wh{Y}(f)\cdot\wh{fX}-2f\cdot\wh{Y}(f)\cdot\wh{X}+2f[\wh{Y},\wh{fX}]
-f^2[\wh{Y},\wh{X}].
\eeas
Hence, similarly as in (\ref{i4}),
\beas\zF([\wh{Y},\wh{f^2X}])&=&
2(\wh{Y}(f)\circ\zc)\cdot\zF(\wh{fX})-2(f\circ\zc)\cdot(\wh{Y}(f)\circ\zc)\cdot\zF(\wh{X})\\
&&+2(f\circ\zc)\cdot\zF([\wh{Y},\wh{fX}])
-(f\circ\zc)^2\cdot\zF([\wh{Y},\wh{X}]).
\eeas
Comparing the above with
$$[\zF(\wh{Y}),\zF(\wh{f^2X})]=[\zF(\wh{Y}),2(f\circ\zc)\cdot\zF(\wh{fX})
-(f\circ\zc)^2\cdot\zF(\wh{X})],
$$
we get easily
\be\label{i7a}
\left(\wh{Y}(f)\circ\zc-\zF(\wh{Y})(f\circ\zc)\right)
\left(\zF(\wh{fX})-(f\circ\zc)\cdot\zF(\wh{X})\right)=0.
\ee
After polarizing with $f:=f+h$ and multiplying both sides by
$\wh{Y}(f)\circ\zc-\zF(\wh{Y})(f\circ\zc)$, we get the identity
\be\label{i8}
\left(\wh{Y}(f)\circ\zc-\zF(\wh{Y})(f\circ\zc)\right)^2
\left(\zF(\wh{hX})-(h\circ\zc)\cdot\zF(\wh{X})\right)=0,
\ee
valid for all $f,h\in\cC(N_1)$ and all $X,Y\in\E_1$. From
(\ref{i8}) we get
$$(\wh{Y}(f)\circ\zc)(q)=(\zF(\wh{Y})(f\circ\zc))(q)$$
for such $q=\zf(p)\in S_2$ for which in no neighbourhood of them
the anchor map is a differential operator of order 0, i.e., for
$q$ which do not belong to
$$S_2^0=\{ \zf(p)\in S_2: \wh{hX}(p')=h(p')\wh{X}(p') \text{\ for all\
}h\in\cC(N_1),\ X\in\E_1\text{\ and\ } p' \text{\ close to \
}p\}.
$$
If, on the other hand, $q\in S_2^0$, then
$\zF(\wh{hX})(q')=(h\circ\zc)(q')\cdot\zF(\wh{X})(q')$ for $q'$
from a neighbourhood of $q$, so that comparing in this
neighbourhood
$$\zF([\wh{Y},\wh{fX}])=(\wh{Y}(f)\circ\zc)\cdot\zF(\wh{X})
+(f\circ\zc)\cdot\zF([\wh{Y},\wh{X}])$$ with
$$[\zF(\wh{Y}),\zF(\wh{fX})]=\zF(\wh{Y})(f\circ\zc)\cdot\zF(\wh{X})
+(f\circ\zc)\cdot[\zF(\wh{Y}),\zF(\wh{X})]$$ we get
$$(\wh{Y}(f)\circ\zc)(q)\cdot\zF(\wh{X})(q)=
\zF(\wh{Y})(f\circ\zc)(q)\cdot\zF(\wh{X})(q),$$ thus
$$(\wh{Y}(f)\circ\zc)(q)=\zF(\wh{Y})(f\circ\zc)(q)$$
also for $q\in S_2^0$.

\ep

\bigskip\noindent
{\bf Remark.}

\smallskip\noindent
(a) For Jacobi-Kirillov bundles with all leaves of the
characteristic foliation (i.e., orbits of $\wh{\E}$) of dimension
$>1$ there is much simpler argument showing that $\zc$ is smooth
than the one using Lemma \ref{l2}. The difficulty in the general
case comes from singularities of the `bivector field' part of the
anchor map and forced us to use Lemma \ref{l2}.

\medskip\noindent
(b) Theorem \ref{t5} has been proven for Lie algebroids in
\cite{GG}, so the new (and difficult) here is the case of
Jacobi-Kirillov bundles with non-trivial `bivector part' of the
bracket. A similar result for the Lie algebras of smooth vector
fields preserving a symplectic or a contact form up to a
multiplicative factor has been proven by H.~Omori \cite{O}. These
Lie algebras are the Lie algebras of locally hamiltonian vector
fields for the Jacobi-Kirillov brackets associated with the
symplectic and the contact form, respectively.

\begin{cor}\

\smallskip\noindent
(a) If the Lie algebras $(\cC(N_i),\{\cdot,\cdot\}_{\beta_i})$ of
the Jacobi contact brackets, associated with contact manifolds
$(N_i,\beta_i)$, $i=1,2$, of class $\cC$, are isomorphic, then the
manifolds $N_1$ and $N_2$ are $\cC$-diffeomorphic.

\medskip\noindent
(b) If the Lie algebras associated with nowhere-vanishing Poisson
structures of class $\cC$ on class $\cC$ manifolds $N_i$, $i=1,2$,
are isomorphic, then the manifolds $N_1$ and $N_2$ are
$\cC$-diffeomorphic.
\end{cor}

\noindent Janusz GRABOWSKI\\Institute of Mathematics\\Polish
Academy of Sciences\\\'Sniadeckich 8\\P.O. Box 21\\00-956 Warsaw,
Poland\\Email: jagrab@impan.gov.pl\\\\
\end{document}